\documentclass[a4paper, 11pt, twoside]{article}
\usepackage[nogytt]{mfpaperstuff2}
\usepackage{tikz-3dplot}
\usetikzlibrary{shapes}
\begin{document}
\newcommand\kst[2]{\operatorname{K}_{#1#2}}

\title{A note on Kostka numbers}
\date{}

\maketitle

\begin{abstract}
We give an elementary proof of a well-known result on Kostka numbers, following a question from Mark Wildon on MathOverflow \cite{mo}. Namely, we show that given partitions $\la,\mu,\nu$ of $n$ with $\mu\dom\nu$, we have $\kst\la\nu\gs\kst\la\mu$.
\end{abstract}

\section{Introduction}

Recall that a \emph{composition} of $n$ is a sequence $\la=(\la_1,\la_2,\dots)$ of non-negative integers which sum to $n$. Given compositions $\la$ and $\mu$ of $n$, we say that $\la$ \emph{dominates} $\mu$ (written $\la\dom\mu$) if $\la_1+\dots+\la_r\gs\mu_1+\dots+\mu_r$ for every $r$.

A composition is a \emph{partition} if it is weakly decreasing. The \emph{Young diagram} of a partition $\la$ is the set
\[
[\la]=\lset{(r,c)\in\bbn^2}{c\ls\la_r},
\]
which we draw as an array of boxes with the English convention (so that $r$ increases down the page, and $c$ from left to right). A \emph{$\la$}-tableau is a function from $[\la]$ to $\bbn$, and we depict a tableau $T$ by drawing $[\la]$ and filling each box with its image under $T$. The \emph{type} of $T$ is the composition $\mu$, where $\mu_i$ is the number of $i$s appearing in the diagram.

A $\la$-tableau is \emph{semistandard} if the entries weakly increase from left to right along rows, and strictly increase down the columns. Given a partition $\la$ of $n$ and a composition $\mu$ of $n$, the \emph{Kostka number} $\kst\la\mu$ is the number of different $\la$-tableaux of type $\mu$.

This note concerns the following well-known result.

\begin{thm}\label{main}
Suppose $\la$ and $\mu$ are partitions of $n$. Then $\kst\la\mu>0$ \iff $\la\dom\mu$.
\end{thm}

The `only if' part of \cref{main} is easy to see: if $T$ is a semistandard $\la$-tableau of type $\mu$, then all the numbers less than or equal to $r$ in $T$ must occur in the first $r$ rows, so $\la_1+\dots+\la_r\gs\mu_1+\dots+\mu_r$. The converse is trickier to prove combinatorially, though a construction is given by the author in \cite{mo}. The objective here is to give an elementary proof of the following result.

\begin{propn}\label{domk}
Suppose $\la,\mu,\nu$ are partitions of $n$ with $\mu\dom\nu$. Then $\kst\la\mu\ls\kst\la\nu$.
\end{propn}

Since obviously $\kst\la\la=1$, this proves the `if' part of \cref{main}. We remark in passing that \cref{domk} (and our proof) works when $\la$ is a skew Young diagram.

\section{The proof of \cref{domk}}

First we require an elementary lemma. Given non-negative integers $x_1,\dots,x_r,a$, let $S(x_1,\dots,x_r;a)$ be the number of ways choosing integers $y_1,\dots,y_r$ such that $0\ls y_i\ls x_i$ for each $i$ and $y_1+\dots+y_r=a$. Now we have the following.

\begin{lemma}\label{elem}
Suppose $x_1,\dots,x_r,a,b$ are non-negative integers, and let $m=x_1+\dots+x_r$. If $|a-\frac m2|\gs|b-\frac m2|$, then $S(x_1,\dots,x_r;a)\ls S(x_1,\dots,x_r;b)$.
\end{lemma}

\begin{pf}
Note first that $S(x_1,\dots,x_r;a)=S(x_1,\dots,x_r;m-a)$, since we have a bijection defined by $y_i\mapsto x_i-y_i$.  So (replacing $a$ with $m-a$ if necessary, and similarly for $b$) we can assume $a\ls b\ls\frac m2$.  Assuming $r\gs 1$ and $x_1\gs1$, we write
\[
S(x_1,\dots,x_r;a)=T(x_1,\dots,x_r;a)+U(x_1,\dots,x_r;a),
\]
where $T(x_1,\dots,x_r;a)$ is the number of ways of choosing the $y_i$ with $y_1=x_1$, and $U(x_1,\dots,x_r;a)$ is the number of ways of choosing the $y_i$ with $y_1<x_1$.  Obviously we have
\[
T(x_1,\dots,x_r;a)=S(x_2,\dots,x_r;a-x_1),\qquad U(x_1,\dots,x_r;a)=S(x_1-1,x_2,\dots,x_r;a)
\]
so it suffices to show that
\[
S(x_2,\dots,x_r;a-x_1)\ls S(x_2,\dots,x_r;b-x_1),\qquad S(x_1-1,x_2,\dots,x_r;a)\ls S(x_1-1,x_2,\dots,x_r;b).
\]
The first of these follows by induction, since $b-x_1$ is at least as close to $(m-x_1)/2$ as $a-x_1$ is.  And the second also follows, since $b$ is at least as close to $(m-1)/2$ as $a$ is.  So we can use induction on $m$.
\end{pf}

Using this, we can prove the following result which is the main ingredient in the proof of \cref{domk}.

\begin{lemma}\label{adjl}
Suppose $i\in\bbn$, $\la$ is a partition of $n$, and $\mu$ is a composition of $n$ with $\mu_i>\mu_{i+1}$. Define a composition $\nu$ by
\[
\nu_i=\mu_i-1,\quad\nu_{i+1}=\mu_{i+1}+1,\quad\nu_j=\mu_j\ \text{ for all other }j.
\]
Then $\kst\la\mu\ls\kst\la\nu$.
\end{lemma}

\begin{pf}
We define an equivalence relation $\sim$ on semistandard $\la$-tableaux by setting $S\sim T$ if all the entries different from $i$ and $i+1$ are the same in $S$ as they are in $T$. We show that within any one equivalence class there are at least as many semistandard tableaux of type $\nu$ as of type $\mu$.

So fix an equivalence class $\calc$, and consider how to construct semistandard tableaux in $\calc$. The positions of the entries different from $i$ and $i+1$ are determined, and we may as well assume there are $\mu_j$ entries equal to $j$ for each $j\neq i,i+1$ (otherwise $\calc$ contains no tableaux of type $\mu$ or $\nu$). We are left with some positions in which to put $i$s and ($i+1$)s -- call these \emph{available} positions. There are at most two available positions in each column, and if there are two, then these must be filled with $i$ and $i+1$. So we need only consider columns having exactly one available position. Given $j\gs1$, let $x_j$ be the number of columns having an available position in row $j$ only; these columns are consecutive, and can be filled in any way with $i$s and $(i+1)$s as long as the $i$ are to the left of the $(i+1)$s, to produce a semistandard tableau.

So choosing a semistandard tableau in $\calc$ amounts to choosing integers $y_1,y_2,\dots$ such that $0\ls x_j\ls y_j$ for each $j$: $y_j$ is just the number of $i$s placed in available positions in row $j$. In order for this semistandard tableau to have type $\mu$, we must have $y_1+y_2+\dots=a$, where $a=\frac12(\mu_i-\mu_{i+1}+x_1+x_2+\dots)$. Similarly, to obtain a semistandard tableau of type $\nu$ we must have $y_1+y_2+\dots=b$, where $b=\frac12(\mu_i-\mu_{i+1}-2+x_1+x_2+\dots)$. Since $\mu_i>\mu_{i+1}$, $b$ is at least as close to $\frac12(x_1+x_2+\dots)$ as $a$ is, so by \cref{elem} there are at least as many tableaux of type $\nu$ in $\calc$ as there are of type $\mu$.
\end{pf}

In order to use \cref{adjl} we need to describe the covers in the dominance order on partitions. We leave the proof of the following results as an easy exercise.

\begin{propn}\label{domcover}
Suppose $\mu$ and $\nu$ are partitions of $n$ with $\mu\doms\nu$. Then $\mu$ covers $\mu$ in the dominance order on partitions (i.e.\ there is no \emph{partition} $\xi$ with $\mu\doms\xi\doms\nu$) \iff one of the following occurs:
\begin{itemize}
\item
for some $i\in\bbn$ we have
\[
\nu_i=\mu_i-1,\quad\nu_{i+1}=\mu_{i+1}+1,\quad\nu_j=\mu_j\ \text{ for all other }j;
\]
\item
for some $i,j\in\bbn$ with $i<j$ we have
\[
\mu_{i+1}=\dots=\mu_j=\mu_i-1,\quad\nu_i=\mu_i-1,\quad\nu_j=\mu_j+1,\quad\nu_k=\mu_k\ \text{ for all other }k.
\]
\end{itemize}
\end{propn}

Informally, $\mu$ covers $\nu$ \iff $\nu$ is obtained by moving one box down and to the right, either to an adjacent row or to an adjacent column.

\begin{pf}[Proof of \cref{domk}]
We may assume $\mu$ covers $\nu$ in the dominance order, and consider the two cases in \cref{domcover}. In the first case it is immediate from \cref{adjl} that $\kst\la\mu\ls\kst\la\nu$. In the second case, define compositions $\xi^{i+1},\dots,\xi^{j-1}$ by
\[
\xi^k_i=\mu_i-1,\quad\xi^k_k=\mu_k+1,\quad,\xi^k_l=\mu_l\ \text{ for all other }l.
\]
Then by \cref{adjl} we have
\[
\kst\la\mu\ls\kst\la{\xi^{i+1}}\ls\cdots\ls\kst\la{\xi^{j-1}}\ls\kst\la\nu.\qedhere
\]
\end{pf}


\begin{thebibliography}{99}

\backrefparscanfalse

\bibitem[MO]{mo}
M.~Wildon, \textit{Is there a short proof that the Kostka number $\kst\la\mu$ is non-zero whenever $\la$ dominates~$\mu$?}, mathoverflow.net/questions/226537.

\end{thebibliography}
\end{document}